\documentclass[11pt,A4paper]{article}

\usepackage[left=32mm,right=32mm,top=30mm,bottom=32mm]{geometry}
\usepackage{labelfig}
\usepackage{epsfig}
\usepackage{epstopdf}

\usepackage{xfrac}

\usepackage[percent]{overpic}

\usepackage{color}
\usepackage{amsthm,amsmath,amssymb}
\usepackage{booktabs}
\usepackage{mathpazo}
\usepackage{microtype}
\usepackage{overpic}
\usepackage{bm}
\usepackage{sectsty}
\usepackage[
	pdftitle={PDFTitle},
	pdfauthor={Hugo Parlier},
	ocgcolorlinks,
	linkcolor=linkred,
	citecolor=linkred,
	urlcolor=linkblue]
{hyperref}

\definecolor{linkred}{rgb}{0.48,0.1,0.05}
\definecolor{linkblue}{RGB}{16, 78, 139}

\usepackage[hang,flushmargin]{footmisc}
\usepackage{enumitem}
\usepackage{titlesec}
	\titlespacing{\section}{0pt}{12pt}{0pt}
	\titlespacing{\subsection}{0pt}{6pt}{0pt}
	
\titlelabel{\thetitle.\quad}

\makeatletter 

\long\def\@footnotetext#1{%
\H@@footnotetext{%
\ifHy@nesting 
\hyper@@anchor{\@currentHref}{#1}%
\else 
\Hy@raisedlink{\hyper@@anchor{\@currentHref}{\relax}}#1%
\fi 
}}

\def\@footnotemark{%
\leavevmode 
\ifhmode\edef\@x@sf{\the\spacefactor}\nobreak\fi 
\H@refstepcounter{Hfootnote}%
\hyper@makecurrent{Hfootnote}%
\hyper@linkstart{link}{\@currentHref}%
\@makefnmark 
\hyper@linkend 
\ifhmode\spacefactor\@x@sf\fi 
\relax 
}%

\ifFN@multiplefootnote%
\renewcommand*\@footnotemark{%
\leavevmode 
\ifhmode 
\edef\@x@sf{\the\spacefactor}%
\FN@mf@check 
\nobreak 
\fi 
\H@refstepcounter{Hfootnote}%
\hyper@makecurrent{Hfootnote}%
\hyper@linkstart{link}{\@currentHref}%
\@makefnmark 
\hyper@linkend 
\ifFN@pp@towrite 
\FN@pp@writetemp 
\FN@pp@towritefalse 
\fi 
\FN@mf@prepare 
\ifhmode\spacefactor\@x@sf\fi 
\relax%
}%
\fi 

\makeatother 

\theoremstyle{plain}
\newtheorem{theorem}{Theorem}[section]

\newtheorem{lemma}[theorem]{Lemma}
\newtheorem{corollary}[theorem]{Corollary}

\theoremstyle{definition}

\newcommand{\floor}[1]{\left\lfloor #1 \right\rfloor}

\newcommand{\R}{{\mathbb R}}

\newcommand{\Hyp}{{\mathbb H}}

\newcommand{\M}{{\mathcal M}}

\newcommand{\Area}{{\rm Area}}
\newcommand{\arcsinh}{{\,\rm arcsinh}}
\newcommand{\arccosh}{{\,\rm arccosh}}

\newcommand{\injrad}{{\rm injrad}}

\sectionfont{\large \sc}
\subsectionfont{\normalsize}

\long\def\symbolfootnote[#1]#2{\begingroup%
\def\thefootnote{\fnsymbol{footnote}}\footnote[#1]{#2}\endgroup}

\def\blfootnote{\xdef\@thefnmark{}\@footnotetext}

\begin{document}

{\Large \bfseries \sc Chromatic numbers of hyperbolic surfaces}

{\bfseries Hugo Parlier\symbolfootnote[1]{\normalsize Research supported by Swiss National Science Foundation grant number PP00P2\textunderscore 153024}, Camille Petit \symbolfootnote[2]{\normalsize Research supported by Swiss National Science Foundation grant number 200021\textunderscore 153599
}
}

{\em Abstract.} 
This article is about chromatic numbers of hyperbolic surfaces. For a metric space, the $d$-chromatic number is the minimum number of colors needed to color the points of the space so that any two points at distance $d$ are of a different color. We prove upper bounds on the $d$-chromatic number of any hyperbolic surface which only depend on $d$. In another direction, we investigate chromatic numbers of closed genus $g$ surfaces and find upper bounds that only depend on $g$ (and not on $d$). For both problems, we construct families of examples that show that our bounds are meaningful. 
\vspace{1cm}

\section{Introduction} \label{s:introduction}

The chromatic number of a graph is the minimum number of colors needed to color its vertices so that any two adjacent vertices are colored differently. Given a metric space $(X, \delta)$ and a real number $d>0$, one can associate a (possibly infinite) graph where vertices are points of $X$ and two vertices are joined by an edge if they are exactly at distance $d$. The chromatic number of this graph gives rise to a notion of chromatic number $\chi((X,\delta),d)$ for a metric space. The particular case when $X$ is the Euclidean plane (and $d=1$ although any choice of $d$ is equivalent) has attracted a particular amount of attention and is called the Hadwiger-Nelson problem (see  \cite{Soifer11, ShelahSoifer1, ShelahSoifer2} and references therein). By exhibiting an explicit coloring coming from a hexagonal tiling, it is not particularly difficult to prove that it is at most $7$. A lower bound of $4$ can be obtained by exhibiting a four-chromatic unit distance graph in the plane, for instance the Moser spindle. Going beyond these two rather elementary bounds is completely open and will probably require either great perseverance or an inspired idea. 

Other metric spaces have been investigated including $n$-dimensional Euclidean space (see for instance \cite{SoiferBook,Taha-Kahle}) and more recently the hyperbolic plane $\Hyp$ \cite{Kloeckner}. Unlike Euclidean spaces, the hyperbolic plane is not invariant by homothety, so a priori the chromatic number depends on a choice of $d$. Bounds for $\chi(\Hyp,d)$ in function of $d$ have been established by Kloeckner, but it is not known whether or not there exists a uniform (independent of $d$)  upper bound. A theorem of de Bruijn-Erd\"os \cite{DeBruijnErdos} says that any given infinite graph can be colored by $k$ colors if and only if all of its finite subgraphs can as well. So, showing that the chromatic number of $\Hyp$ can be made arbitrarily large amounts to exhibiting ``subgraphs" of $\Hyp$ (by which we mean geometric copies of finite graphs) with arbitrarily large chromatic number. For the moment, the best known lower bound is only $4$. 

Our main focus is on the more general setup of hyperbolic surfaces (not necessarily the plane). We ask only that they be complete hyperbolic surfaces. Our first results is an upper bound on the chromatic number which only depends on $d$. 
\begin{theorem}\label{thm:A1}
There exists a constant $C_1>0$ such that for every number $d>0$ every complete hyperbolic surface $S$ satisfies 
$$
\chi(S,d) \leq C_1\, e^d.
$$
\end{theorem}
In first instance, our upper bound, exponential in $d$, seems particularly weak in comparison with the linear upper bound for the hyperbolic plane. But in fact we exhibit, for any $d>0$, surfaces with a chromatic number that is at least exponential in $\sfrac{d}{2}$. 
\begin{theorem}\label{thm:A2}
There exists a constant $C_2>0$ and a family of complete hyperbolic surfaces $S_d$, $d>0$ so that
$$
\chi(S_d,d)\geq C_2\, e^{\sfrac{d}{2}}.
$$
\end{theorem}
The optimal constants have growth that are exponential in $\alpha d$ for some $\alpha$ in between $\frac{1}{2}$ and $1$. Determining the exact value for $\alpha$ could be an interesting problem. The construction and proofs are quite elementary and only require some basic tools about hyperbolic geometry and trigonometry. 

Given a metric space $(X,\delta)$, there is a natural way of associating a chromatic number that doesn't depend on $d$. One defines the chromatic number of $(X,\delta)$ as
$$
\chi((X,\delta)) := \sup \{ \chi((X,\delta), d) \,:\, d>0\}.
$$
When $(X,\delta)$ is the Euclidean plane, this is simply the chromatic number discussed previously. When $(X,\delta)$ is the hyperbolic plane, we've seen that it is unknown whether this quantity is finite. In the particular case when $(X,\delta)$ is a compact Riemannian manifold however, then by a compactness argument, this quantity is always finite. Now if one has  a natural family of compact manifolds, one can study the supremum of this quantity over the whole family. As an example, consider the following problem (which we don't know the answer to): among all $2$ dimensional flat tori, which one has the largest chromatic number? By the theorem of de Bruijn and Erd\"os mentioned previously, this quantity is an upper bound for the chromatic number of the plane - in fact the chromatic number of {\it any} torus is an upper bound (this will be explained in the preliminaries). 

For closed hyperbolic surfaces, one can ask the same question. In this case, for each genus, we get a moduli space of isometry types of surfaces. Our second set of results are about bounds on the chromatic numbers that only depend on the genus and not on the individual geometries. We begin with our upper bounds. 

\begin{theorem}\label{thm:B1}
There exists a constant $C_3>0$ such that for every integer $g\geq 2$ every closed hyperbolic surface $S$ of genus $g$ satisfies
$$
\chi(S)\leq C_3\, g.
$$
\end{theorem}

Again, one could ask whether there is not a universal upper bound (which doesn't depend on genus) but we exhibit families of surfaces that provide the following lower bounds. 
\begin{theorem}\label{thm:B2}
There exists $C_4>0$ and a family of closed hyperbolic surfaces $S_g$, where $S_g$ has genus $g$, such that
$$\chi(S_g)\geq C_4\, \sqrt g.
$$
\end{theorem}
So again, we show that there exists a constant $\beta$, lying somewhere in between $\frac{1}{2}$ and $1$ such that the optimal upper bound on chromatic numbers behaves like $g^{\alpha}$. Whether there is relation between the constant $\alpha$ described above and this constant $\beta$ remains to be seen. Interestingly, our results rely on the celebrated result of Ringel and Youngs about the genus of complete graphs (which already provided an answer to another graph coloring problem, namely the Heawood conjecture). We also mention that one could ask the same questions, and apply some of the same techniques, to other moduli spaces, such as hyperbolic surfaces with punctures but for the sake of clarity, we've restricted our study to closed surfaces. 

One by-product of our lower bounds in the above theorem is an example of hyperbolic surface which has infinite chromatic number. It is not literally a corollary of the theorem but can be directly obtained using the same building blocks as the ones we need in the proof of Theorem \ref{thm:B2}.

\begin{corollary}\label{cor:infinite}
There exists a hyperbolic surface $Z$ such that $\chi(Z)=\infty.$ 
\end{corollary}

This article is organized as follows. After a preliminary section in which we introduce tools on the geometry of hyperbolic surfaces, we have two main sections. The first of these is about bounds of $d$-chromatic numbers in terms of $d$ and the second about closed surfaces and bounds on chromatic numbers in terms of the genus. 

{\bf Acknowledgement.} Both authors would like to thank Bill Balloon for inspiration.

\section{Preliminaries}

We'll use this preliminary section to introduce definitions and notations, and also to give a short description of how we're thinking about hyperbolic surfaces and some of the properties we'll use throughout the paper. 

\subsection{Chromatic numbers}

For a metric space $(X,\delta)$ we define its chromatic number $\chi((X,\delta),d)$ relative to a distance $d>0$ (or simply $d$-chromatic number) as the minimal number of colors needed to color all points of $X$ such that any $x,y \in X$ with $\delta(x,y) = d$ are colored differently. We call a $d$-coloring a coloring of $(X,\delta)$ where any $x,y \in X$ with $\delta(x,y) = d$ are colored differently. When the metric space consists in the vertices of a graph, distance to edge distance and $d=1$, this is the usual definition of the chromatic number of the graph. 

Equivalently, one can define the chromatic number of a metric space using the usual chromatic number of graphs by associating a graph to the metric space as follows. Given a metric space $(X, \delta)$ and a real number $d>0$, we construct a graph whose vertices are points of $X$ and we place an edge between points  if they are exactly at distance $d$. 

For certain metric spaces, the choice of $d$ is crucial; for others, such as $n$-dimensional Euclidean space, any choice of $d$ is equivalent. This prompts the following definition. We define the chromatic number $\chi((X,\delta))$ of a metric space $(X,\delta)$ to be the quantity
$$
\chi((X,\delta)):= \sup \{ \chi((X,\delta),d) \, :\, d>0\}.
$$
As examples, for the Euclidean and hyperbolic planes, the following inequalities are known:
$$4 \leq \chi(\R^2) \leq 7 \mbox{ and }4 \leq \chi(\Hyp) \leq \infty.$$

The theorem of de Bruijn and Erd\"os \cite{DeBruijnErdos} , mentioned in the introduction, is the following:
\begin{theorem} Any infinite graph $G$ can be colored by $k$ colors if and only if all of its finite subgraphs can be colored by $k$ colors. 
\end{theorem}

As an example of an application of this theorem, consider any Euclidean flat torus $T$ of dimension $n$. Then we claim that
$$
\chi(T) \geq \chi(\R^n).
$$
To see this, apply the de Bruijn-Erd\"os theorem to $\R^n$. There is thus a finite set of points of $\R^n$ that realize $\chi(\R^n)$. They can be made to lie in a ball of arbitrarily small size. For a fixed $T$, one can find any Euclidean ball of sufficiently small size that is isometrically embedded inside $T$. Thus the chromatic number of $T$ is at least the chromatic number of $\R^n$. Whether or not these quantities are equal - or finding an upper bound on $\chi(T)$ in function of $\chi(\R^n)$ - could be an interesting problem. In particular, understanding the behavior of
$$
\sup \{ \chi(T) \, :\, T \mbox{ is a flat $n$-dimensional torus}\}
$$
in function of $n$. 

Problems for $n$-dimensional tori can often be translated to analogous questions to hyperbolic surfaces of genus $g\geq 2$ where the genus plays the part of the dimension. These sets of are two natural generalizations of the set of $2$-dimensional flat tori. We denote $\M_g$ the {\it moduli space} of genus $g\geq 2$ hyperbolic surfaces which can be thought of the space of closed genus $g\geq 2$ hyperbolic surfaces up to isometry. From the chromatic number $\chi(S)$ of a hyperbolic surface $S$ of genus $g\geq 2$, one can study
$$
\sup \{ \chi(S) \, :\, S \in \M_g\}
$$
as a function of $g$. Investigating this quantity is one of the main goals of the article. In the next section we begin by properly defining which types of hyperbolic surfaces we're interested in and some of the tools we'll need in the sequel. 

\subsection{Hyperbolic surfaces and their thick-thin decomposition}

A hyperbolic surface is a surface locally isometric to the hyperbolic plane $\Hyp$. A surface is said to be complete if it geodesically complete as a Riemannian manifold. We will generally be concerned with complete hyperbolic surfaces but to construct them we will sometimes use surfaces with (simple) geodesic boundary. We denote $d_S$ the distance function for a surface $S$. 

The simplest such surface is a funnel: topologically an infinite half-cylinder. Such a surface can be obtained as follows: we quotient $\Hyp$ by a hyperbolic element to obtain an infinite cylinder. Such a cylinder has a unique simple closed geodesic (the quotient of the axis of the hyperbolic element by its action) which separates the cylinder into two half cylinders. One of these half cylinders is the funnel we're talking about.

These funnels provide a way of going from a surface with simple geodesic boundary to a complete surface by pasting funnels of the appropriate boundary length to the surface with boundary. Note that the original surface with boundary is a convex subset of the full surface and is sometimes referred to as the convex core. From a dynamical point of view, everything interesting on the surface happens within the convex core.  

On hyperbolic surfaces, there is a unique geodesic representative in every isotopy class of a simple closed curve (which does not surround a disk or a cusp). Simple closed geodesics of hyperbolic surface have an associated {\it collar} which is a tubular neighbourhood around the geodesics which can be described as follows. 
 
\begin{lemma}[Collar lemma]
 Let $\gamma$ be a simple closed geodesic on a complete hyperbolic surface $S$. Then the set 
 $$
 {\mathcal C}(\gamma) := \{ x \in S \, : \, d_S(x, \gamma) \leq w(\gamma) \}
 $$
 where
 $$
 w(\gamma)=\arcsinh\left( \frac{1}{\sinh\left(\frac{\ell(\gamma)}{2}\right)}\right)
 $$
 is an embedded cylinder isometric to $[-w(\gamma), w(\gamma)] \times {\mathcal S}^1$ with the Riemannian metric
 $$
 ds^2 = d\rho^2 + \ell^2(\gamma) \cosh^2(\rho) dt^2.
 $$
 \end{lemma}
This version of the collar lemma, and many of the basic facts we state, can be found in \cite{BuserBook}. 

Furthermore, if any two simple closed geodesics are disjoint, then so are their collars. Immediate consequences of the collar lemma include the fact that any two simple closed geodesics of length less than $2 \arcsinh(1)$ are always disjoint. One place collars naturally appear is when a surface is decomposed into its {\it thick} and {\it thin} parts. For a given $\varepsilon>0$, we can separate a surface $S$ into its $\varepsilon$-thick part, namely
$$
\widehat S^{\varepsilon} := \left\{ x \in S \,:\,  \injrad(x) > \varepsilon \right\}
$$
where $\injrad(x)$ is the injectivity radius of $S$ at the point $x$ and its $\varepsilon$-thin part
$$
S\setminus \widehat S^{\varepsilon} := \left\{ x \in S \,:\,  \injrad(x) \leq \varepsilon \right\}.
$$
Because any two simple closed geodesics of length less than $2\arcsinh(1)$ cannot intersect, when $\varepsilon < \arcsinh(1)$ the set $S\setminus \widehat S^{\varepsilon} $ (if it is not empty!) consists of a collection of cylinders. These cylinders are either collars of a certain width around a simple closed geodesic or are neighbourhoods of cusps. 

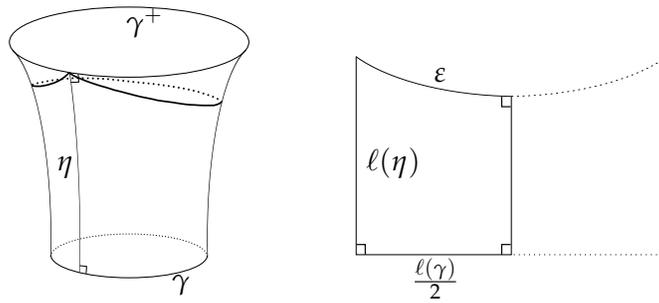
\begin{figure}[h]
\leavevmode \SetLabels
\L(.25*.4) $\eta$\\
\L(.35*.0) $\gamma$\\
\L(.52*.4) $\ell(\eta)$\\
\L(.56*.0) $\frac{\ell(\gamma)}{2}$\\
\L(.58*.7) $\varepsilon$\\
\L(.31*.87) $\gamma^+$\\
\endSetLabels
\begin{center}
\AffixLabels{\centerline{\epsfig{file =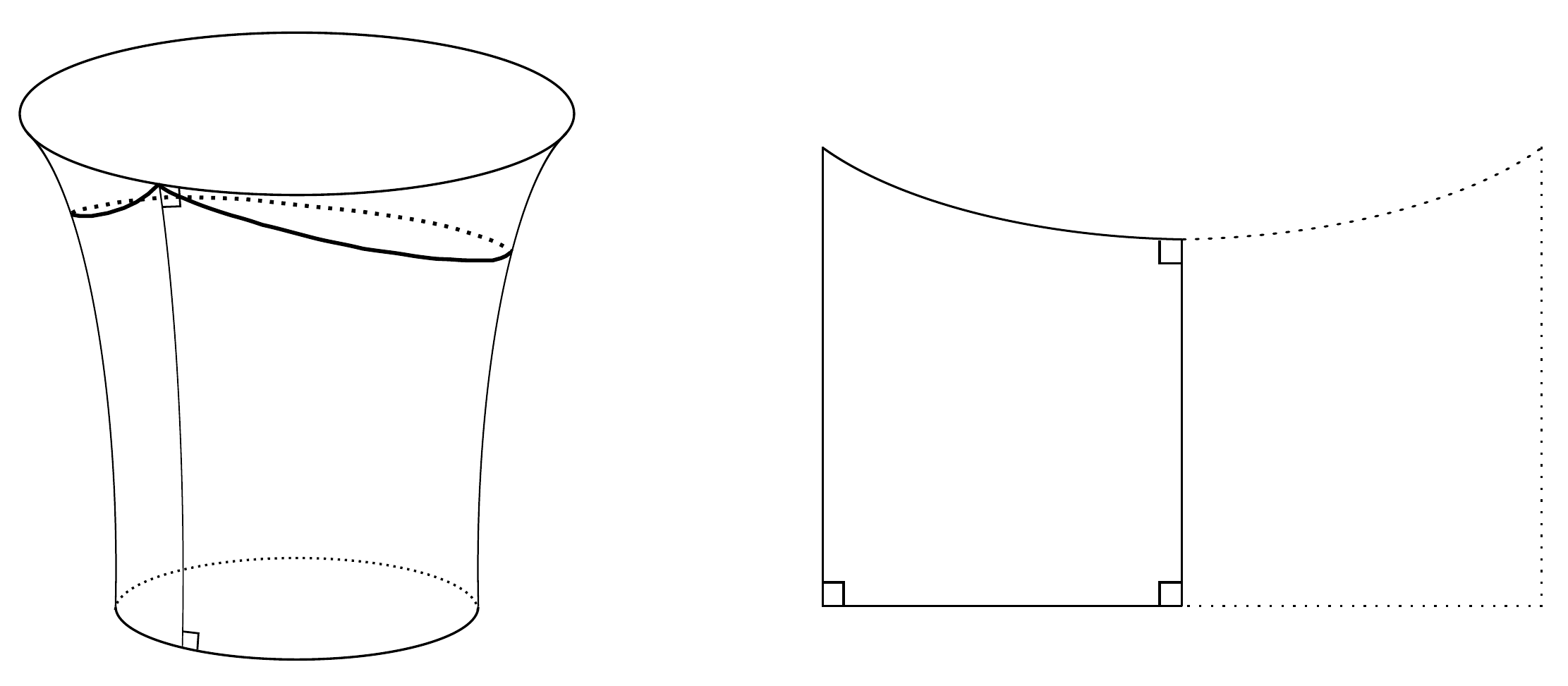,width=9.0cm,angle=0} }}
\vspace{-30pt}
\end{center}
\caption{Computing $K_C$} \label{fig:KC}
\end{figure}

Consider one of them, say $C$, that contains a simple closed geodesic $\gamma$. Any point on one of its boundary curves is the base point of an embedded geodesic loop of length $2 \varepsilon$. From this we can deduce that the two boundary curves, say $\gamma^+$ and $\gamma^-$, of $C$ are smooth curves, both of equal length, and both parallel lines to the unique closed geodesic in their homotopy class $\gamma$. To see this, observe that as these loops are all of equal length, and all parallel to $\gamma$, the angle formed by any of these loops at the base point must always be the same. This is simply because the length of $\gamma$ can be computed as a function of this angle and the length of the loop (which is always $2 \varepsilon$). Similarly, the distance in between any of the boundary points and $\gamma$ is always equal. As such, for any $C$, there exists a $K_C$ which only depends on $C$ (or alternatively $\varepsilon$ and $\ell(\gamma)$) such that 
$$
\partial C = \{ x \in S \,:\, d_S(x, \gamma) = K_C\}=\gamma^+ \cup \gamma^-.
$$
Let us compute the value of $K_C$ in function of $\ell(\gamma)$ and $\varepsilon$. Fix a point on the boundary of $C$ and consider a distance path $\eta$ to $\gamma$. Cutting along the loop of length $2\varepsilon$, $\gamma$ and $\eta$ gives hyperbolic quadrilateral such as Figure \ref{fig:KC}. 
This quadrilateral can be divided into two quadrilaterals with three right angles as in the figure with we can compute. By standard hyperbolic trigonometry we have 
$$
\sinh(\varepsilon) = \sinh\left(\frac{\ell(\gamma)}{2}\right) \cosh(\ell(h))
$$
from which we obtain 
$$
K_C= \ell(\eta) = \arccosh\left( \frac{\sinh(\varepsilon) }{\sinh\left(\frac{\ell(\gamma)}{2}\right) }\right).
$$
It's interesting to compare this value to the width of the collar from the collar lemma. In particular note that the difference
\begin{equation}\label{eqn:collardifference}
w(\gamma) - K_C= \arcsinh\left( \frac{1}{\sinh\left(\frac{\ell(\gamma)}{2}\right)}\right) -  \arccosh\left( \frac{\sinh(\varepsilon) }{\sinh\left(\frac{\ell(\gamma)}{2}\right) }\right)
\end{equation}
is a positive number because $\varepsilon < \arcsinh(1)$. Furthermore it reaches its minimum in $\ell(\gamma)=0$. When $\varepsilon = \arcsinh(1)$, this minimum is exactly $0$. 

The geodesic $\gamma$ also divides $C$ into two parts, $C^{+}, C^{-}$, whose other boundary curves are $\gamma^{+}$ and $\gamma^{-}$. 

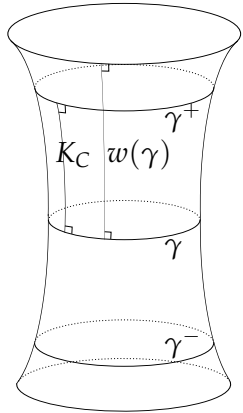
\begin{figure}[h]
\leavevmode \SetLabels
\L(.44*.6) $K_C$\\
\L(.485*.6) $w(\gamma)$\\
\L(.535*.68) $\gamma^+$\\
\L(.535*.18) $\gamma^-$\\
\L(.535*.4) $\gamma$\\
\endSetLabels
\begin{center}
\AffixLabels{\centerline{\epsfig{file =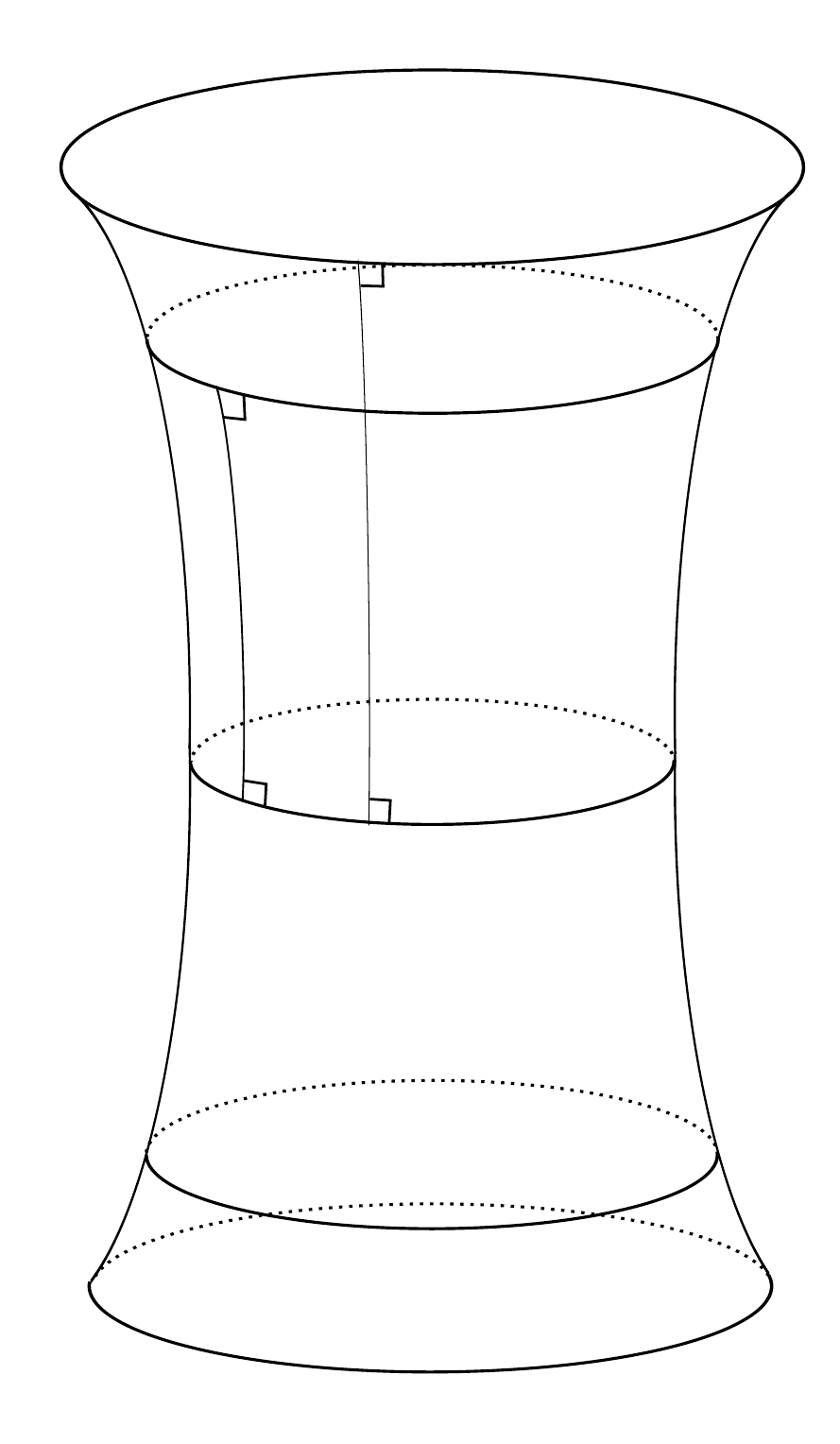,width=3.5cm,angle=0} }}
\vspace{-30pt}
\end{center}
\caption{The collar around a simple closed geodesic $\gamma$} \label{fig:}
\end{figure}

We now put a further restriction on $\varepsilon$ which will be useful in what follows. Although $C$ may not be convex, provided $\varepsilon>0$ is small enough, both $C^{+}, C^{-}$ will be. We suppose now that 
$$
\varepsilon \leq \arcsinh\left(\frac{1}{\sqrt{2}}\right)
$$
and we'll see that this condition ensures convexity. By computing the limit when $\ell(\gamma)=0$ in the quantity in Equation \eqref{eqn:collardifference}, we get a lower bound on this quantity which is equal to $\frac{\log(2)}{2}$. What will be crucial in the argument that follows is that twice this value is less than $\arcsinh\left(\frac{1}{\sqrt{2}}\right)$.

We now turn our attention to the convexity of $C^{+}$. To see this, we begin by observing that the convex hull of either of the boundary curves of $C$ lies in its respective half: for example for any $x,y\in \gamma^{+}$, the shortest geodesic between $x,y$ lies in $C^{+}$. The shortest path inside $C^{+}$ between $x$ and $y$ is of length strictly less than half $\ell(\gamma^{+})$. So if this is not the minimal geodesic path, then there is another shorter path between them. 

\begin{figure}[h]
\leavevmode \SetLabels
\L(.08*.53) $x$\\
\L(.2*.54) $y$\\
\L(.13*.52) $\gamma^+$\\
\L(.455*.68) $x$\\
\L(.52*.51) $y$\\
\L(.53*.67) $\gamma^+$\\
\L(.55*.425) $\gamma$\\
\L(.825*.465) $x$\\
\L(.88*.25) $y$\\
\L(.87*.455) $\gamma^+$\\
\L(.92*.2) $\gamma$\\
\endSetLabels
\begin{center}
\AffixLabels{\centerline{\epsfig{file =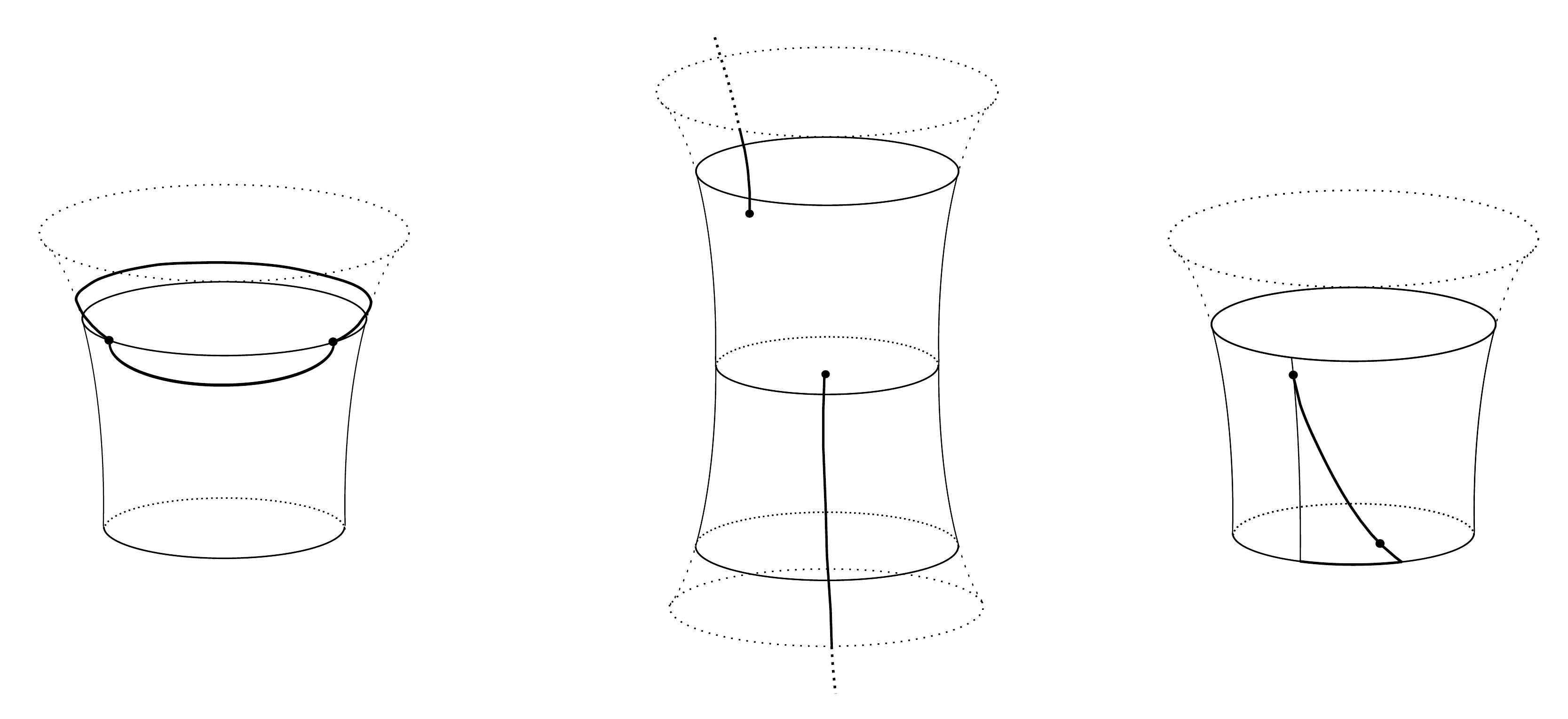,width=15.0cm,angle=0} }}
\vspace{-30pt}
\end{center}
\caption{Three types of potential distance paths} \label{fig:HalfCylinderConvexity}
\end{figure}
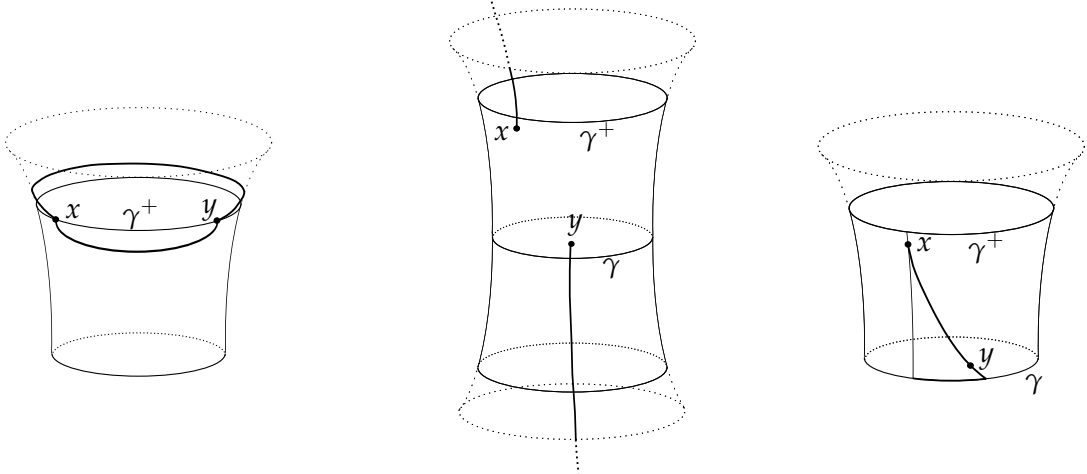

To show this never occurs, we fix a point $x\in\gamma^+$ and move the point $y\in\gamma^+$ away from $x$ until this occurs for the first time. The path obtained lies entirely outside the interior of $C^{+}$ and these two simple geodesic paths form a curve (see the left most configuration of Figure \ref{fig:HalfCylinderConvexity}). On a hyperbolic surface there are no geodesic bigons so the corresponding curve is non-trivial and has a geodesic representative $\tilde{\gamma}$ which in turn has its own collar. By construction the two collars $C(\gamma)$ and $C(\tilde\gamma)$ intersect,  which is impossible by the version of the collar lemma described in the preliminaries. We have reached a contradiction and shown that the convex hull of $\gamma^{+}$ is entirely contained in $C^{+}$. 

Now the only way $C^{+}$ can be non-convex is if there is a distance path which leave $C^{+}$ and returns through $C^{-}$ (see the middle case of Figure \ref{fig:HalfCylinderConvexity}). As such it will have length at least the width of $C^{-}$ which is equal to $K_C$. In addition it will have spent some time in the thick part of $S$: by the estimates given above this will add at least $
\log(2)$ to its length. So it is of length at least $K_C +\log(2)$. The shortest path between the two points  that lies entirely inside $C^{+}$ is of length at most 
 $$K_C +\frac{\ell(\gamma)}{2}\leq K_C  + \arcsinh\left(\frac{1}{\sqrt{2}}\right) < K_C + \log(2)$$
 and this provides a contradiction. 

\section{Bounding the chromatic number in function of $d$}\label{s:distance}

In this section we prove Theorem \ref{thm:A1} from the introduction. We begin by proving a universal upper bound on the chromatic number which only depends on the parameter $d$. Then, for every $d>0$, we exhibit a surface $S_{d}$ which has chromatic number at least $C\, e^{\sfrac{d}{2}}$ for some universal constant $C>0$.

\subsection{Upper bounds}

To prove our upper bounds, we will need to lift any complete hyperbolic surface $S=\Hyp / \Gamma$ to its universal cover $\Hyp$ and then construct a $\Gamma$ invariant coloring of $\Hyp$. 

To begin, for a fixed value $r>0$, we consider a maximal set $\Delta_r\subset S$ of points such that if $x,y \in \Delta_r$ and $x\neq y$ then $d_S(x,y)>r$. By construction this set satisfies the following two properties:
\begin{itemize}[itemsep=2ex,leftmargin=0.5cm]
\item For $x,y\in\Delta_r$, $x\neq y$, $B_{\sfrac{r}{2}}(x) \cap B_{\sfrac{r}{2}}(y) = \emptyset$.
\item $S = \cup_{x \in \Delta_r} B_{r}(x) $
\end{itemize}
Here, $B_{r}(x)$ denotes the closed ball center at $x$ of radius $r$.
Note that we don't ask that the sets $B_{\sfrac{r}{2}}(x)$ be embedded balls in $S$. 

We remark however that for any $\rho$, a set $B_{\rho}(x)$ lifts to a union of embedded balls on $\Hyp$ (not necessarily disjoint - for instance if $\rho$ is larger than the diameter of $S$, then the lift is $\Hyp$). Furthermore, if $B,B'$ are two disjoint balls on $S$, then any of their lifts are also disjoint and the distance $d_S(B,B')$ is simply the minimum of the distances of their lifts in $\Hyp$. This is all essentially in the definitions of a cover or the universal cover but we emphasize it as it will be crucial in what follows. 

Now given $d>0$ we want to color $S$ by balls of radius $r_0$ so we set $r_0:= \min \{ \frac{2d}{5}, \arcsinh (1)\}$. In particular each ball is of diameter strictly less than $d$. We then consider a $\Delta_{r_0}$ as described above.

We now endow the set $\Delta_{r_0}$ with a graph structure $G$ as follows. Vertices are points of $\Delta_{r_0}$ and two vertices $x,y$ share an edge if there exists $x'\in B_{r_0}(x)$ and $y'\in B_{r_0}(y)$ such that $d_S(x',y') = d$. Our strategy is to bound the degree of $G$ by a function of $d$ and $r_0$. 

To do this we lift a point $x\in\Delta_{r_0}$ to the universal cover. We denote $\tilde x\in \pi^{-1}(x)$, where $\pi:\Hyp \to S$ is the covering map. We observe that for any $\rho>0$ the set  $B_\rho(\tilde{x})$ (which lies in $\Hyp$) covers the set $B_\rho(x)$ (and of course belongs to its preimage).

For $x\in\Delta_{r_0}$ we want to bound the degree $\deg(x)$ of $x$ in $G$. We compute an upper bound on this cardinality using $\tilde x$: it is bounded by the number of $\tilde y \in \tilde\Delta_{r_0} := \pi^{-1}(\Delta_{r_0})$ such that there exists $x'\in B_{r_0}(\tilde x)$ and $y'\in B_{r_0}(\tilde y)$ which satisfy
$$d_\Hyp (x',y')=d.$$
The balls of radius $\sfrac{r_0}{2}$ around any such $\tilde y$ are disjoint and must lie entirely in the annulus $A$ centered at $\tilde x$, of inner radius $d-\frac{5}{2}r_0$ and outer radius $d+\frac{5}{2}r_0$. The area of a ball of radius $\rho$ in the hyperbolic plane is
$$
4\pi \sinh^2\left(\frac{\rho}{2}\right)
$$
so we have
$$
4\pi | \{\tilde y \in \tilde\Delta_{r_0}: B_{\sfrac{r_0}{2}}(\tilde y) \subset A\}| \sinh^2\left(\frac{r_0}{4}\right) < 4\pi \left(\sinh^2\left(\frac{d+\frac{5}{2}r_0}{2}\right) - \sinh^2\left(\frac{d-\frac{5}{2}r_0}{2}\right)  \right).
$$
Using this we deduce the following bound on $| \{\tilde y \in \tilde\Delta_{r_0}: B_{\sfrac{r_0}{2}}(\tilde y) \subset A\}|$ which in turn bounds the degree of any point of $G$:
\begin{eqnarray*}
\deg(G) \leq | \{\tilde y \in \tilde\Delta_{r_0}: B_{\sfrac{r_0}{2}}(\tilde y) \subset A\}| &\leq& \frac{\sinh^2\left(\frac{d+\frac{5}{2}r_0}{2}\right) - \sinh^2\left(\frac{d-\frac{5}{2}r_0}{2}\right) }{\sinh^2\left(\frac{r_0}{4}\right)}\\
&=& \frac{\sinh(\frac{5}{2}r_0)}{\sinh^2\left(\frac{r_0}{4}\right)} \sinh(d).
\end{eqnarray*}

Recalling the definition of $r_0$, we obtain $\deg(G)\leq \phi(d)$, where
\begin{equation*}
  \phi(d):=\left\{ \begin{array}{ll}
   \frac{\sinh^2(d)}{\sinh^2(\sfrac{d}{10})}, & d\leq 10 \arcsinh(1)\\
   \sinh(10\arcsinh(1)) \cdot \sinh(d), & d\geq 10 \arcsinh(1).
  \end{array}\right.
 \end{equation*}

We can now deduce our upper bounds. We use Brooks' theorem on graph coloring which asserts that a graph of degree at most $D$ can be colored with at most $D+1$ colors. 

We thus obtain a coloring of $G$  with at most $\phi(d)+1$ colors. It induces a $d$-coloring of $S$ as follows. We color each ball $B_{r_0}(x), x\in\Delta_{r_0}$ with the color corresponding to the vertex $x$ in $G$. If a point belongs to several balls we choose one of the colors of the balls it belongs to arbitrarily as its color. This proves the upper bound $\chi(S,d)\leq \phi(d)+1$, which in particular gives
$$
\chi(S,d)\leq C\, e^d
$$
for some constant $C>0$, as stated in Theorem \ref{thm:A1}.
Note that our $d$-coloring of $S$ lifts to a $\Gamma$ invariant $d$-coloring of $\Hyp$.

\subsection{Lower bounds}

The goal is to give, for any $d>0$, a surface $S_d$ with a $d$-chromatic number that satisfies the lower bound of Theorem \ref{thm:A2}.

Before going to the general construction, we begin by constructing a family of surfaces $S_{d_N}$ for a discrete set of $d_N$s with $d_N \to \infty$ as $N\to \infty$. The general construction that follows retains many of the key properties of this simpler construction.

For any integer $N\geq 3$, there is a unique ideal regular hyperbolic polygon with $N$ sides. It has a well defined center point and a number of self isometries including rotations of angle $\frac{2\pi}{N}$ around this center point. It also has a unique maximally embedded disk (we'll compute its radius $R_N$ in the sequel), centered in the center point. By symmetry the disk touches each of the $N$ sides in points which we'll call the {\it midpoints} of the sides. 

\begin{figure}[h]
\leavevmode \SetLabels
\L(.51*.47) $v_j$\\
\L(.62*.68) $v_k$\\
\L(.53*.61) $d_N$\\
\endSetLabels
\begin{center}
\AffixLabels{\centerline{\epsfig{file =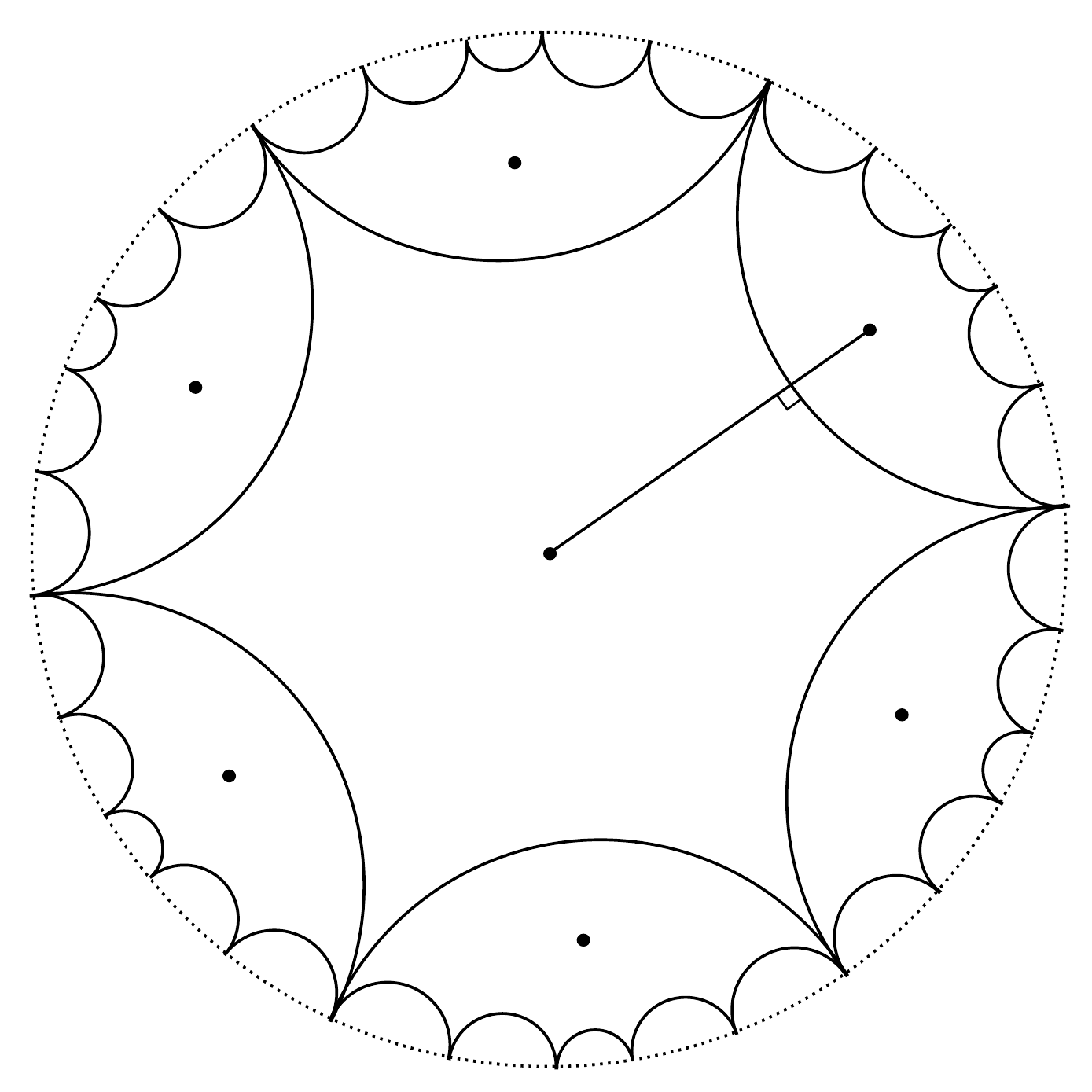,width=6cm,angle=0} }}
\vspace{-30pt}
\end{center}
\caption{Constructing $S_{d_N}$} \label{fig:RegularPolygon}
\end{figure}

For any $N$ we consider $N+1$ copies of this unique ideal regular hyperbolic polygon with $N$ sides and glue the sides in pairs (see Figure \ref{fig:RegularPolygon}). The only thing we ask of this gluing is that it must obey the following two rules: every two distinct polygons share exactly one side and the sides are pasted in their midpoints. The result is a connected finite area surface with at least one cusp. 

We denote $v_1,\hdots,v_N$ the center points of the polygons. These will be the vertices of an embedded complete graph of $N$ vertices. 

Our first claim is that on the resulting surface
$$
d(v_j,v_k) = 2 R_N
$$
for $k\neq j$. This is simply because the distance between a center point $v_k$ and any of the sides of the polygon is $R_N$. So any path between two distinct vertices must pass through one of the sides of the polygons and as such has at least length $2R_N$. Between any two distinct vertices there is a unique path of length $2R_N$ given by the concatenation of the radii and this proves the claim. So we have geometric embedding of $K_N$ with edge length $2 R_N$. 

\begin{figure}[h]
\leavevmode \SetLabels
\L(.255*.1) $v_j$\\
\L(.715*.1) $v_k$\\
\L(.48*.0) $d_N$\\
\L(.355*.13) $\sfrac{\pi}{N}$\\
\L(.595*.13) $\sfrac{\pi}{N}$\\
\endSetLabels
\begin{center}
\AffixLabels{\centerline{\epsfig{file =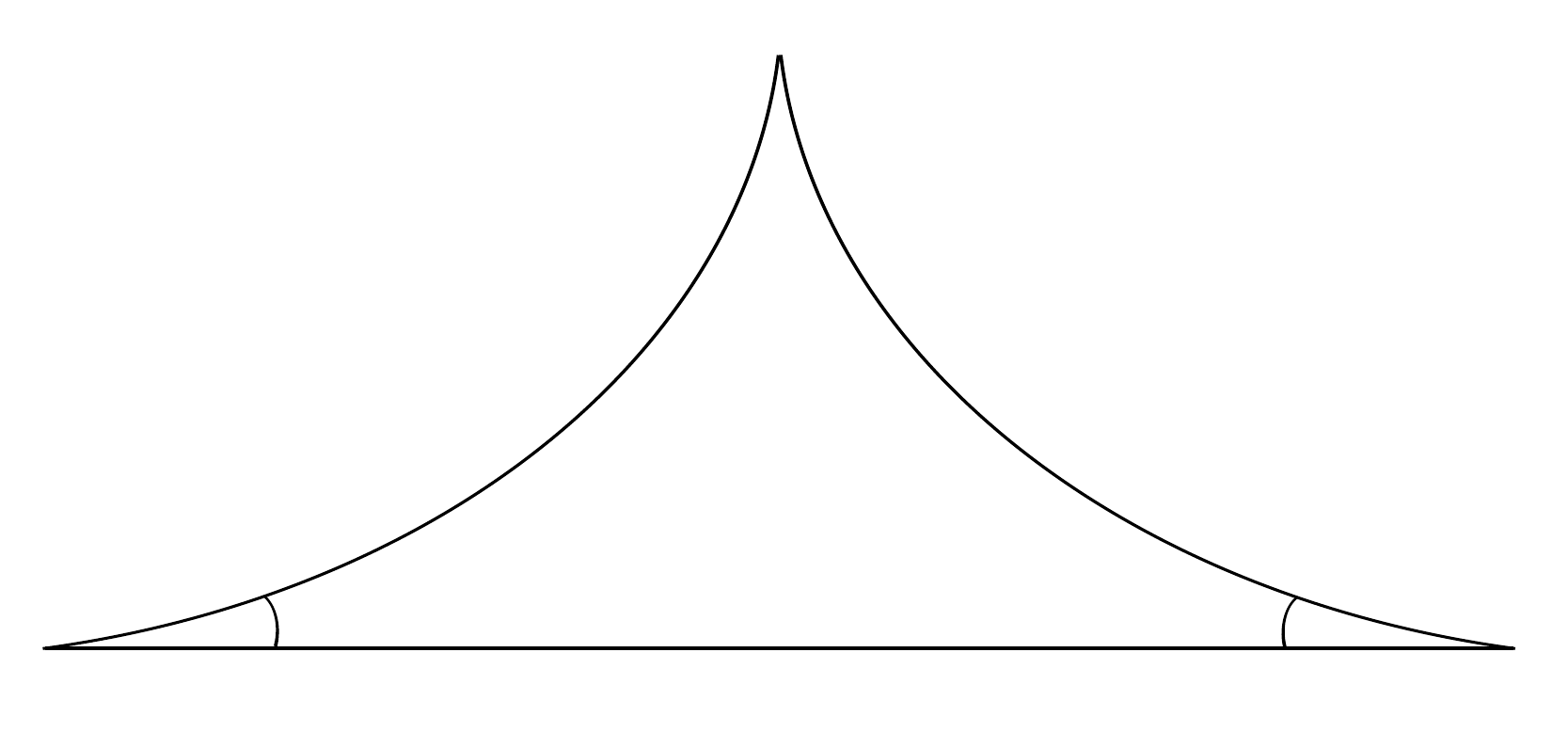,width=7cm,angle=0} }}
\vspace{-30pt}
\end{center}
\caption{Computing $d_N$} \label{fig:IsocelesTriangle}
\end{figure}
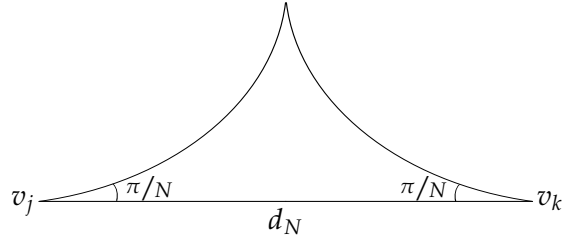

Let us set $d_N:= 2 R_N$ and compute its value. We consider a triangle formed by any two distinct vertices $v_j,v_k$ and an ideal point as in Figure \ref{fig:IsocelesTriangle}.

It has angles $0, \frac{\pi}{N},  \frac{\pi}{N}$ so by standard hyperbolic trigonometry the following holds:

\begin{eqnarray*}
\cosh(d_N) &=& \frac{1+ \cos^2\left(\sfrac{\pi}{N} \right)}{\sin^2\left(\sfrac{\pi}{N} \right)}\\
&=& \frac{2}{\sin^2\left(\sfrac{\pi}{N} \right)} - 1.
\end{eqnarray*}
From this we obtain
$$
d_N= \arccosh\left(  \frac{2}{\sin^2\left(\sfrac{\pi}{N} \right)} - 1\right).
$$
We observe that $d_N$ grows asymptotically like $2 \log(N)$ as $N$ goes to infinity. Since from our construction we obtain an embedded complete graph with $N$ vertices formed by the centers $\{ v_1,\cdots,v_{N+1}\}$ and with edges length $d_N$, we get the lower bound $\chi(S_{d_N},d_N)\geq N$,
which in turn gives
$$
\chi(S_{d_N},d_N)\geq C \cdot e^\frac{d_N}{2}
$$
for some constant $C>0$.
So this example provides the correct lower bound but only works for a discrete set of values $d_N$. 

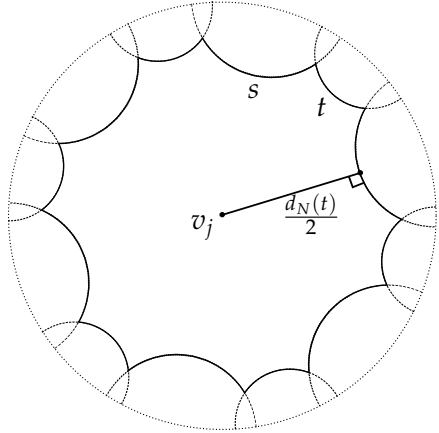
\begin{figure}[h]
\leavevmode \SetLabels
\L(.47*.47) $v_j$\\
\L(.58*.72) $t$\\
\L(.52*.76) $s$\\
\L(.55*.48) $\frac{d_N(t)}{2}$\\
\endSetLabels
\begin{center}
\AffixLabels{\centerline{\epsfig{file =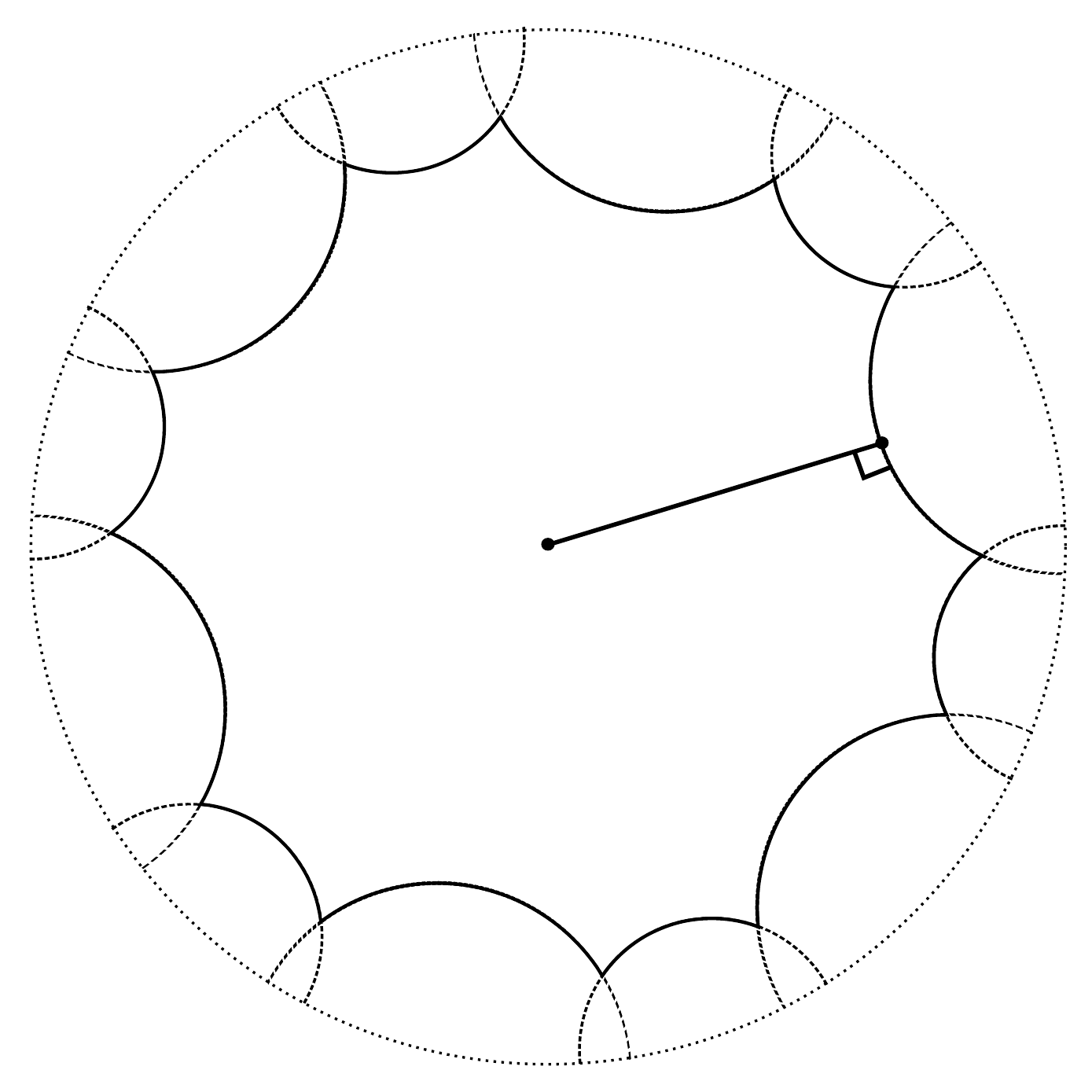,width=6cm,angle=0} }}
\vspace{-30pt}
\end{center}
\caption{A semi-regular right angled $2N$-gon} \label{fig:SemiRegularPolygon}
\end{figure}

We now adapt this construction to construct a surface for every $d\geq d_3$. To do this we replace the ideal polygons in the above construction with semi-regular right angled $2N$-gons with every second side of length $t$ for some $t > 0$; we ask that they have a rotational symmetry of order $N$ which permutes the sides of length $t$ (and thus the $N$ remaining sides as well), as in Figure \ref{fig:SemiRegularPolygon}. These $N$ remaining sides will be of some length $s$ which only depends on $t$ (for fixed $N$). We can see the ideal polygon as the limit case of these polygons  when $t \to 0$. These polygons again have a natural center given by the center of the largest embedded ball inside the polygon. 

We perform the same construction as above but the sides of length of $t$ play the part of the ideal points. More precisely we take $N$ copies of the above polygon and we arbitrarily paste the polygons along their $N$ sides of length $s$ where the only rule is that the resulting surface is orientable and any two distinct polygons are pasted along a single side. Here we don't have to worry about the ``shear" parameter as we are pasting two equal segments together. For each $N$ we get a family of surfaces (with parameter $t$) which has boundary curves. We add hyperbolic funnels to the boundary curves to get complete hyperbolic surfaces $M_N(t)$.

As before we obtain an embedded complete graph with $N$ vertices formed by the centers of the polygons. The distances between these centers now depends on the parameter $t$. For the same reason as in the ideal case, the unique distance paths between the centers is the concatenation of the radial distance paths from the centers to the sides of each polygon pasted together. We denote this distance $d_N(t)$. 

We only need two facts about $d_N(t)$: first of all  $\lim_{t\to 0^{+}} d_N(t) = d_N$ where $d_N$ is as defined above; secondly $d_N(t)$ is a continuous (monotonous) function satisfying $d_N(t) \to \infty$ as $t \to \infty$. The first fact is by construction and the second is a direct consequence of the fact that if $t$ becomes arbitrarily large, $s$ becomes arbitrarily close to $0$ so $d_N(t)$ becomes arbitrarily large. This can be seen more explicitly by hyperbolic trigonometry by relating $t$ and $d_N(t)$:
$$
\cosh\left(\frac{t}{2}\right)=\cosh\left(\frac{d_N(t)}{2}\right) \sin\left(\frac{\pi}{N}\right).
$$

We can now explain how we can associate one of these examples to any $d>0$. We begin by choosing the smallest integer $N$ such that 
$$
d < d_{N+1}.
$$
In particular either $N=2$ (in which the problem is trivial) or 
$$
d_N \leq d < d_{N+1}
$$
for some $N\geq 3$.
We've already constructed the examples for $d=d_N$ so we suppose that $d>d_N$. The  properties of the functions $d_N(t)$ explained above imply that there exists a $t_d>0$ such that
$$
d_N(t_d) = d.
$$
We set the example surface to be $S_d:=M_N(t_d)$ which has $d$-chromatic number bounded below by $N$ by construction. From our previous computations for the ideal surfaces we know
$$
d < d_{N+1} =  \arccosh\left(  \frac{2}{\sin^2\left(\sfrac{\pi}{N+1} \right)} - 1\right).
$$
from which we can deduce the lower bound of Theorem \ref{thm:A2}, that is,
$$
\chi(S_d,d) \geq C\, e^\frac{d}{2}
$$
for some constant $C>0$.

\section{Bounding the chromatic number in function of the genus}\label{s:genus}

In this section we prove Theorems \ref{thm:B1} and \ref{thm:B2} from the introduction. We begin by proving a universal upper bound on the chromatic number which only depends on the genus. We then construct a family of surfaces to prove the lower bound. 

\subsection{Upper bounds}

The first idea for the upper bound is surprisingly simple and is close to what did in a previous section. Let $S$ be a surface of genus $g$. 

We begin by fixing a constant $r_0>0$ which we will specify later but which satisfies $d> 2r_0$. 

Recall the $\frac{r_0}{2}$-thick part  $\widehat S^{\sfrac{r_0}{2}}$of $S$ is defined by:
$$
\widehat S^{\sfrac{r_0}{2}} := \left\{ x \in S \,:\, \injrad(x) > \frac{r_0}{2} \right\}.
$$
We will choose $r_0$ such that $S \setminus \widehat S^{\sfrac{r_0}{2}}$ is a collection of cylinders. On $\widehat S^{\sfrac{r_0}{2}}$ we consider $\Delta_{r_0}$ a maximal set of points on $\widehat S^{\sfrac{r_0}{2}}$ such that if $x,y \in \Delta_{r_0}$ and $x\neq y$ then $d_S(x,y) > r_0$. As the set is maximal we have
$$
\widehat S^{\sfrac{r_0}{2}} \subset \cup_{x\in \Delta_{r_0}} B_{r_0}(x).
$$
and for $x,y\in \Delta_{r_0}$, $x\neq y$, $B_{\sfrac{r_0}{2}}(x) \cap B_{\sfrac{r_0}{2}}(y) = \emptyset$. This allows us to bound the number of points in $\Delta_{r_0}$. We have
$$
\Area\left( \cup_{x\in \Delta_{r_0}} B_{r_0}(x) \right) < \Area(\widehat S^{\sfrac{r_0}{2}}) \leq \Area(S) = 4\pi (g-1)
$$
from which, using the formula for the area of a ball in $\Hyp$, we deduce
$$
4\pi | \Delta_{r_0} | \sinh^2\left( \frac{r_0}{4} \right) < 4 \pi (g-1)
$$
and thus
$$
 | \Delta_{r_0} | < \frac{g-1}{\sinh^2\left( \frac{r_0}{4} \right) }.
$$
We color points of $\widehat S^{\sfrac{r_0}{2}}$ by giving each ball centered in a point of $\Delta_{r_0}$ and of radius $r_0$ a different color. If a point belongs to several balls we choose one of the colors of the balls it belongs to arbitrarily as its color. As balls are of diameter $<d$, it provides a $d$-coloring of $\widehat S^{\sfrac{r_0}{2}}$.

What remains to be colored are points lying in the cylinders comprising $S\setminus \widehat S^{\sfrac{r_0}{2}}$. Consider a cylinder $C$ that lies in this set. To color $C$ we will divide it into sections of diameter less than $d$ and color the sections. 

\begin{figure}[h]
\leavevmode \SetLabels
\L(.155*.23) $h$\\
\L(.27*.84) $\gamma^+$\\
\L(.3*.1) $d'$\\
\L(.54*.37) $h$\\
\L(.69*.79) $\frac{\ell(\gamma^+)}{2}$\\
\L(.7*.32) $d'$\\
\endSetLabels
\begin{center}
\AffixLabels{\centerline{\epsfig{file =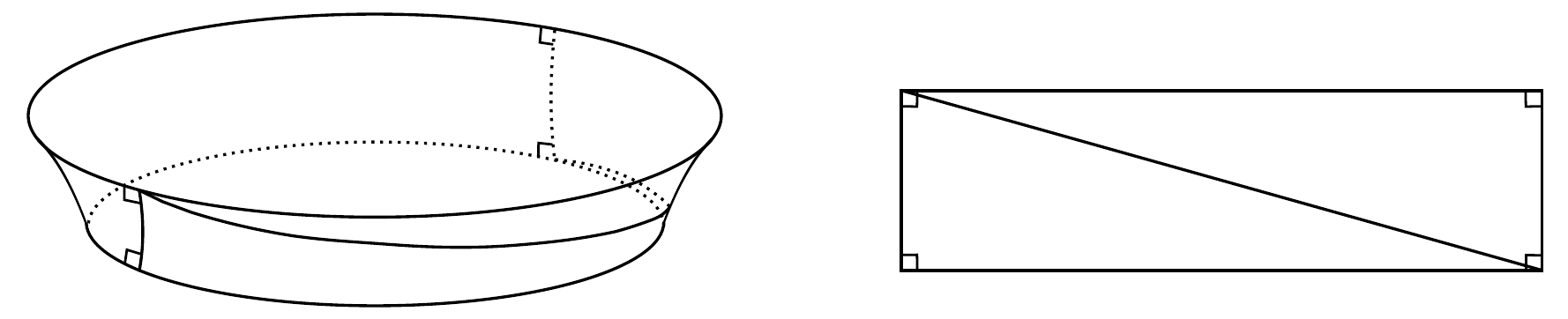,width=12.0cm,angle=0} }}
\vspace{-30pt}
\end{center}
\caption{A section and its diameter $d'$} \label{fig:Section}
\end{figure}
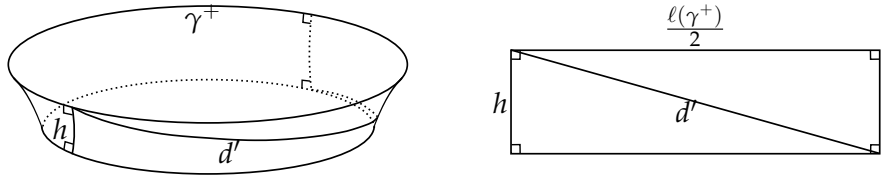

Any point on one of its boundary curves is the base point of an embedded geodesic loop of length $r_0$. As seen in the preliminaries, the two boundary curves $\gamma^+$ and $\gamma^-$ of $C$ are smooth curves, both of equal length, and both parallel lines to the unique closed geodesic $\gamma$ in their homotopy class. There also exists a constant $K_C$ which only depends on $C$ such that
$$
\partial C = \{ x \in S \,:\, d_S(x, \gamma) = K_C\}.
$$
The geodesic $\gamma$ divides $C$ into two convex subsets $C^{+}, C^{-}$. Using the convexity of the half-collars, we can define our sections. Sections are {\it slices} of the half-collars delimited by lines parallel to $\gamma$ in the following way. We want each section to be of diameter less than $d$ but very close to $d$. We define the {\it height} of each section to be the distance between the boundary curves. 
Note that the diameter of a section is realized by pairs of points, one on each of the boundaries of the section, similarly to diametrically opposite points on a Euclidean cylinder (see Figure \ref{fig:Section}). 

\begin{figure}[h]
\leavevmode \SetLabels
\L(.5*-.01) $\gamma$\\
\L(.6*.92) $\gamma^+$\\
\L(.455*.645) $h$\\
\endSetLabels
\begin{center}
\AffixLabels{\centerline{\epsfig{file =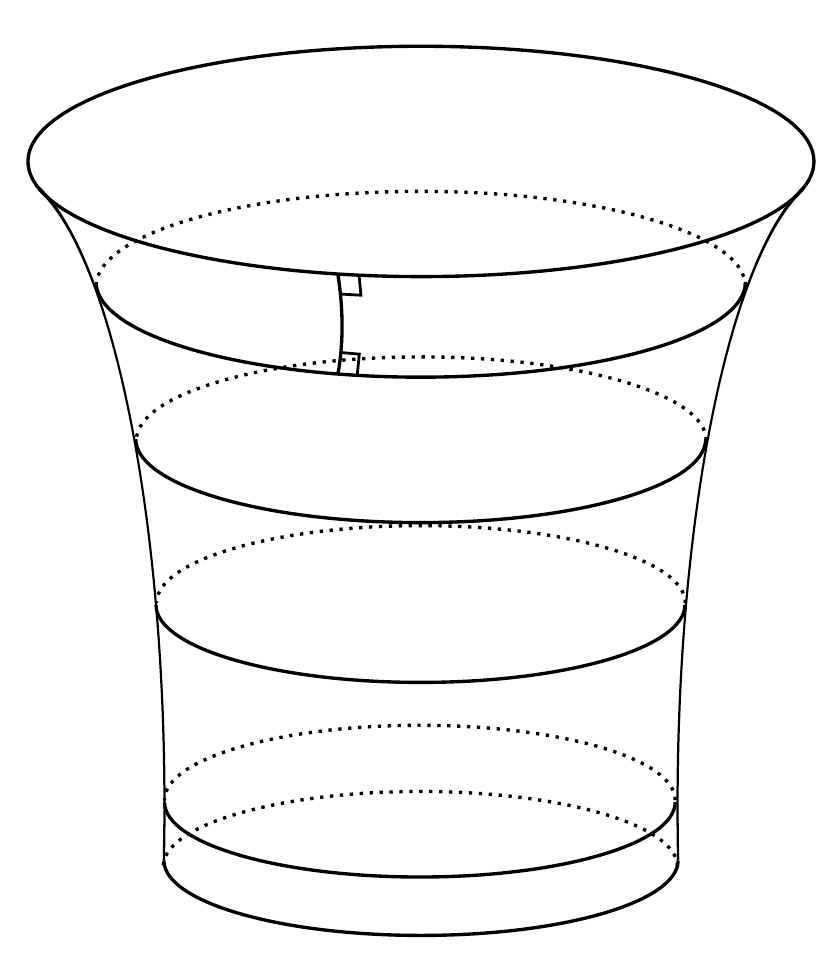,width=4.0cm,angle=0} }}
\vspace{-30pt}
\end{center}
\caption{Slicing a half cylinder into sections} \label{fig:Slices}
\end{figure}

For $C^{+}$ (we proceed analogously for $C^{-}$) the first section we construct is the one with $\gamma^{+}$ as the top boundary curve. Denote $d'$ its diameter with $d-d'$ (arbitrarily) small. There is a triangle with sides of length $h$, $\frac{\ell(\gamma^+)}{2}$ and $d'$ as in Figure \ref{fig:Section}. Thus the height $h$ satisfies the following inequality:
$$
h + \frac{\ell(\gamma^{+})}{2}> d'
$$
From this and the fact that $\ell(\gamma^{+}) \leq r_0$ we have 
$$
h > d' - \frac{r_0}{2}
$$
and thus using $r_0 \leq \frac{d}{2}$ we can choose $d'$ so that 
$$
h > \frac{d}{2}.
$$
We continue to slice $C^{+}$ in sections of diameter close to $d$ iteratively from the top down (see Figure \ref{fig:Slices}). 
Although we've only proved that we can make the height greater than $\frac{d}{2}$ for the first section, to show this was a bound on the ``width" (half the length of a boundary curve). As we move closer and closer to $\gamma$, the boundary curves become smaller and smaller and so the above argument continues to work. Thus we can make the subsequent all of height at least $\frac{d}{2}$. It stops working once we reach $\gamma$ so we don't get a lower bound on the height of the last section (but we won't need one). 

From these sections we create a graph: each section is a vertex and we relate two vertices by an edge if there exist two points, one in each corresponding section, at distance $d$. By the above properties a section is related to at most two subsequent sections and two preceding sections so the graph is of degree at most $4$. From a coloring of the graph we obtain a $d$-coloring of $C^{+}$ by coloring points in each section in the color of the corresponding vertex (boundary points between two sections can be colored by either of the two colors). As the graph is of degree at most $4$, at most $5$ colors are required to $d$-color $C^{+}$. Analogously, we can $d$-color $C^{-}$ with $5$ as well. As we are only interested in the rough growth of the number of colors, although we can clearly $d$-color $C$ with a total of at most $10$ colors. 

Now, there is at most $3g-3$ cylinders so we can $d$-color the totality of the cylinders with at most $10(3g-3)$ colors.

We can now conclude that at most 
$$
 \frac{g-1}{\sinh^2\left( \frac{r_0}{4} \right)}+ 10(3g-3)
 $$
 colors are sufficient to $d$-color any $S$. The upper bound in Theorem \ref{thm:B1} can be deduced by a simple manipulation of the above term setting $r_0 = 4 \arcsinh(1)$ for $d\geq 8 \arcsinh(1)$ and from Theorem \ref{thm:A1} when $d < 8 \arcsinh(1)$.

\subsection{Lower bounds}

The goal is to obtain lower bounds by constructing geometric embeddings of complete graphs with small genus.

We begin with the Ringel and Youngs (\cite{RingelYoungs}) topological embedding of $K_n$ into a surface $M_n$ of genus $g_n$ where
\begin{equation}\label{eqn:RY}
g_n = \floor{ \frac{(n-3)(n-4)}{12}}
\end{equation}
and $g_n$ is the smallest possible genus in which we could embed $K_n$. This embedding $\varphi: K_n \to M_n$ has the following property: for all $n\neq 0$ satisfying $n\equiv 0 \mod 12 $,  
$$
M_n \setminus \varphi(K_n)
$$
is a collection of triangles (see \cite{TerryWelchYoungs}). This topological embedding will serve as a blue print for constructing hyperbolic surfaces with chromatic number roughly root of the genus. As for our lower bounds of Theorem \ref{thm:A2}, we begin by a slightly easier construction before showing how to make it work in the general case; specifically we first construct a family of surfaces with growing genus before constructing a family with a surface in every genus. 

Fix an integer $N$ such that $N+1 \equiv 0 \mod 12$. Our first construction consists in replacing the triangles in $M_{N+1} \setminus \varphi(K_{N+1})$ by equilateral hyperbolic triangles with all three angles equal to $\frac{2\pi}{N}$. These triangles are our building blocks and they are what will change in the more general construction which will follow.

There is a unique such triangle and its three equal side lengths can be directly computed using hyperbolic trigonometry. 

$$
\ell_N = \arccosh \left( \frac{\cos^2(\frac{2\pi}{N}) + \cos(\frac{2\pi}{N})}{\sin^2(\frac{2\pi}{N})} \right).
$$

This construction gives a family of well defined smooth hyperbolic surfaces $S_{g_{N+1}}$ of genus $g_{N+1}$ for a family of $N \to \infty$. What we claim is that $\chi( S_{g_{N+1}}, \ell_{N}) \geq N+1$. This will follow from the geometric embedding of $K_{N+1}$ with edge distance $\ell_N$. 

\begin{figure}[h]
\leavevmode \SetLabels
\L(.475*.7) $\frac{\ell_N}{2}$\\
\endSetLabels
\begin{center}
\AffixLabels{\centerline{\epsfig{file =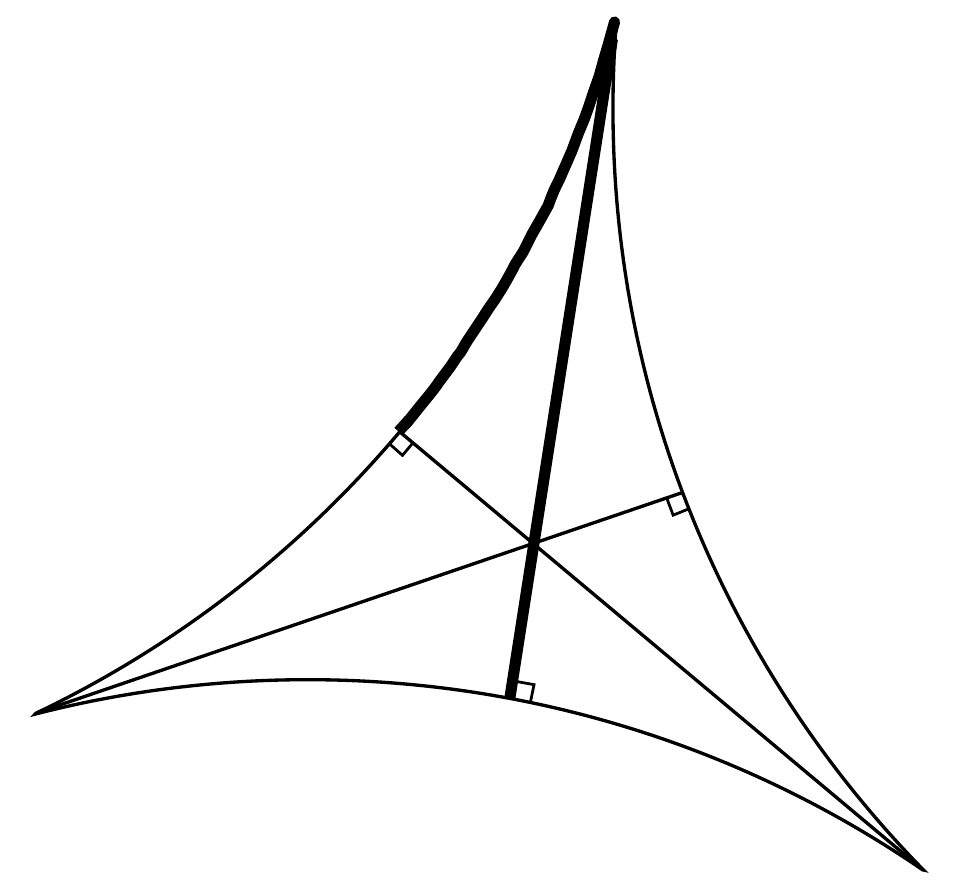,width=5.0cm,angle=0} }}
\vspace{-30pt}
\end{center}
\caption{An equilateral triangle} \label{fig:Triangle}
\end{figure}

To show that the embedding is indeed geometric we need to show that the images of the vertices of $K_{N+1}$ are all at pairwise distance at least $\ell_N$ (by construction they are at distance at most $\ell_N$). Consider a simple geodesic path between two distinct vertices of $\varphi(K_{N+1})$ that is not the side of one of the triangles. We orient this path and look at it as a concatenation of simple segments that each pass through individual triangles. What we claim is that the first of the segments (and hence by symmetry the last) has length strictly greater than $\frac{\ell_N}{2}$. To see this we look at the geometry of the individual triangle. The segment leaves from a vertex and so, as it is geodesic, must leave the triangle through the opposite side. It's length is then at least the minimal distance between a vertex and the opposite side. This distance, in {\it any} equilateral hyperbolic triangle, is strictly greater than half of the length of one of the sides (see Figure \ref{fig:Triangle}).

From this, the length of any geodesic path between distinct vertices which is not the side of a triangle is strictly greater than $\ell_N$. So the embedding is geometric and from Equation \eqref{eqn:RY} with $n=N+1$ we get 
$$
\chi(S_{g_{N+1}},\ell_N) \geq N+1 \geq \sqrt{12g_{N+1} +72}.
$$
So the family of surfaces $S_{g_{N+1}}$ has the lower bound on chromatic number that we are looking for but doesn't have a surface in every genus. We shall need to modify the above construction to get a surface in every genus. 

\begin{figure}[h]
\leavevmode \SetLabels
\L(.4*.58) $\sfrac{\pi}{N}$\\
\L(.57*.38) $\sfrac{\pi}{N}$\\
\L(.51*.95) $\frac{t}{3}$\\
\endSetLabels
\begin{center}
\AffixLabels{\centerline{\epsfig{file =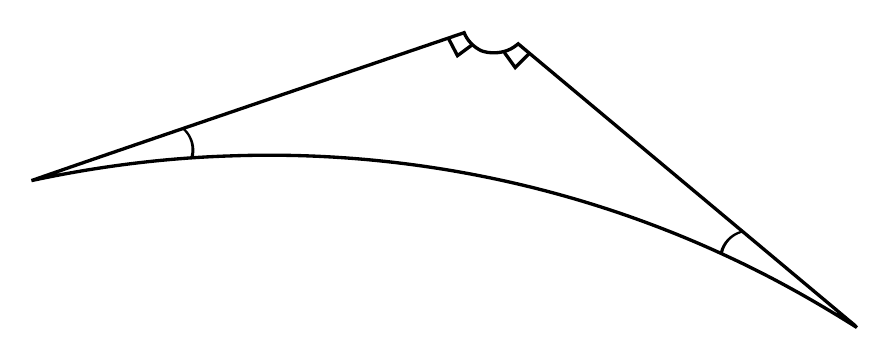,width=6.0cm,angle=0} }}
\vspace{-30pt}
\end{center}
\caption{The quadrilateral used to build the one holed triangle} \label{fig:Quad}
\end{figure}

We begin by changing the building blocks. From equilateral triangles we revert to ``equilateral triangles" with a single interior boundary curve of length $t$. More precisely, for given integer $N$ and given $t>0$ there is a unique ``one holed" triangle with a rotational symmetry of order $3$, one simple geodesic boundary curve of length $t$ and another boundary curve consisting of a triangle with three angles equal to $\frac{2\pi}{N}$. One constructs such a triangle by gluing three copies of a certain quadrilateral: this quadrilateral is the unique quadrilateral with two right angles and a side between them of length $\frac{t}{3}$ and two other angles equal to $\frac{\pi}{N}$ and a central symmetry as is Figure \ref{fig:Quad}.

\begin{figure}[h]
\leavevmode \SetLabels
\endSetLabels
\begin{center}
\AffixLabels{\centerline{\epsfig{file =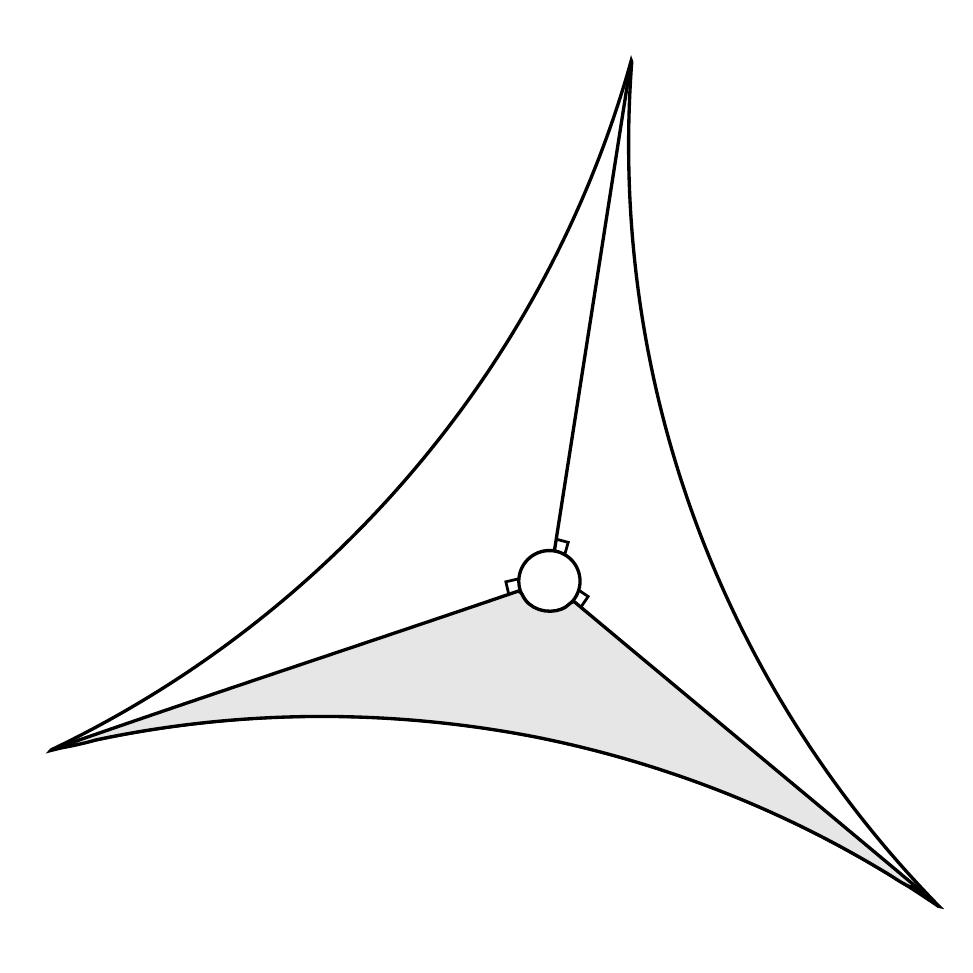,width=5.0cm,angle=0} }}
\vspace{-30pt}
\end{center}
\caption{The one holed triangle} \label{fig:OneHoledTriangle}
\end{figure}

Now by taking three copies of this quadrilateral and pasting them as in Figure \ref{fig:OneHoledTriangle} one obtains the desired building block.

We now analyse the geometry of this one holed triangle in more detail. We begin by fixing $t$ to some small value; exactly which value we choose is of no importance but it will be crucial that it be sufficiently small to satisfy a certain property we shall exhibit in what follows. The fact that we can choose a uniform $t$ independently of $N$ will be important in our infinite genus surface we construct at the very end.

\begin{figure}[h]
\leavevmode \SetLabels
\L(.48*.7) $\frac{\ell}{2}$\\
\L(.53*.55) $a$\\
\L(.48*.41) $c$\\
\endSetLabels
\begin{center}
\AffixLabels{\centerline{\epsfig{file =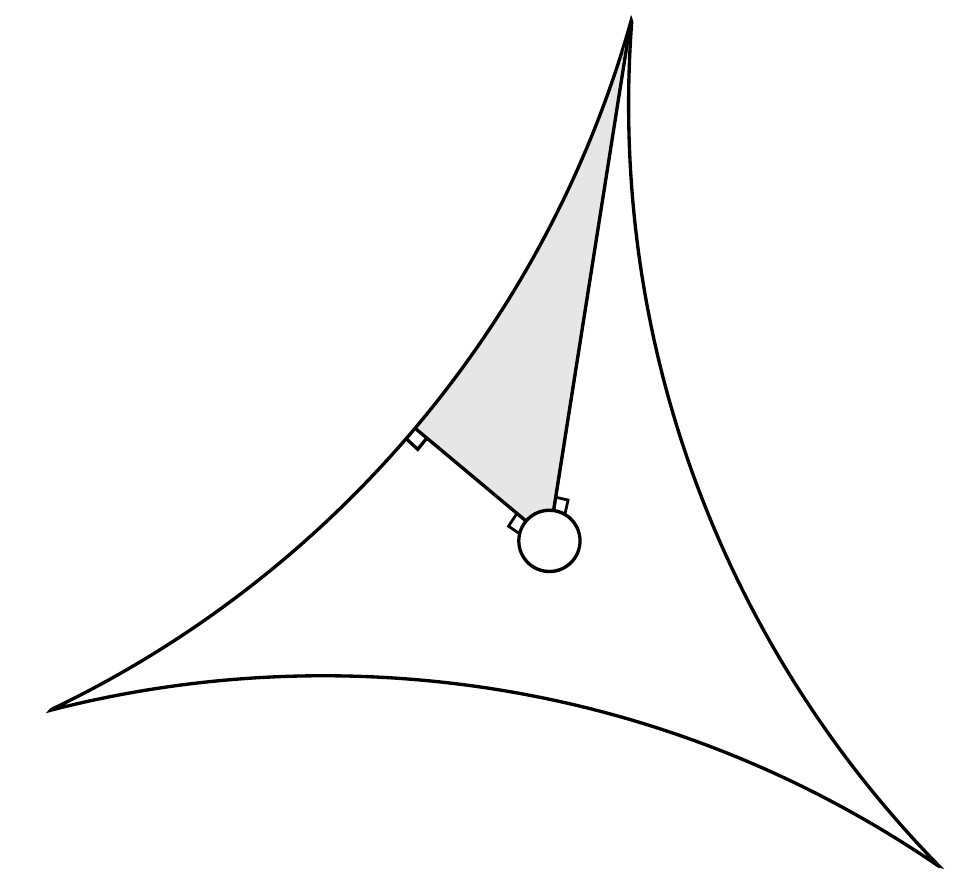,width=5.0cm,angle=0} }}
\vspace{-30pt}
\end{center}
\caption{$a$, $\frac{\ell}{2}$ and $c$} \label{fig:SmallQuad}
\end{figure}

We begin by looking, in the one holed triangle, at the distance $a$ between any of the three vertices and the center hole. We're interested in how this distance relates to the side length, which we denote $\ell$. In particular, we claim that, provided $t$ is small enough, then $a > \frac{\ell}{2}$. Indeed, by looking the quadrilateral highlighted in Figure \ref{fig:SmallQuad}, the values $a,\ell$ and $t$ satisfy the following equality:
$$
\sinh\left(\frac{\ell}{2}\right) = \sinh\left(\frac{t}{6}\right) \cosh(a).
$$

From this we can deduce that 
$$
a = \arccosh\left( \frac{\sinh(\frac{\ell}{2})}{\sinh(\frac{t}{6})}\right)
$$
Since $N\geq 11$, we can guarantee that $\ell$ has a certain length (indeed $\ell$ is certainly longer than $\ell_N$, the side length of the corresponding equilateral triangle with the same angles but {\it no} holes). Now for any $t$ such that $\sinh\left(\frac{t}{6}\right) < 1$, the quantity
$$
 \arccosh\left( \frac{\sinh(\frac{\ell}{2})}{\sinh(\frac{t}{6})}\right) - \frac{\ell}{2}
$$
is positive for large enough $\ell$. Fix for instance $t>0$ such that 
$$
\sinh\left(\frac{t}{6}\right) = \frac{1}{4},
$$
which garanties that $a - \frac{\ell}{2} >0$. 

We now consider a simple geodesic path between a vertex and any other side of the triangle. We claim such a path is of length at least $\frac{\ell}{2}$: indeed such a path must pass through the segment marked $c$ in Figure \ref{fig:SmallQuad}. As such it is of length at least $\frac{\ell}{2}$. 

For fixed $N$, again using hyperbolic trigonometry, we can compute the length $\ell'_N$ of the side of our one holed triangle:
$$
\ell'_N = 2 \arccosh \left( \frac{\cosh(\frac{t}{6})}{\sin(\frac{\pi}{N})} \right).
$$
We now return to the global construction. Again, we paste copies of the one holed triangle using the blueprint provided by the embedding of $K_{N+1}$ in $M_{N+1}$. We obtain a surface with boundary curves, one for each triangle, and all of length $t$. We shall complete the surface by gluing something on each of the boundary curves, but let us remark already that in whatever fashion we do this, the properties of our building blocks imply that the resulting surface has a geometric embedding of $K_{N+1}$ with edge lengths $\ell'_N$. Indeed, any simple path between vertices that is not the side of a triangle can be decomposed into simple segments, each of which lies on a single one holed triangle. What we claim is that, if the path is oriented, the first and last segments are both of length greater than $\frac{\ell}{2}$. Indeed it is either a geodesic path between the vertex and a side of the triangle - so is of length strictly greater than $\frac{\ell}{2}$ as shown above - or it is a path between the vertex and the hole, again of length strictly greater than $\frac{\ell}{2}$. We can conclude that the embedding is geometric. 

For given $N$, this construction gives us a surface of genus $g_{N+1}$ with boundary curves with the property that no matter how we complete it, the resulting surface has a geometrically embedded copy of $K_{N+1}$ for an appropriate edge length. We denote these surfaces $F_N$. 

Let us describe what type of surfaces we can build by pasting this surface in different ways. We shall use it in two ways, the first of which is to construct closed surfaces of different genus to fill the gaps left in the earlier construction. 

{\it \underline{Closing the gaps}}

The number of boundary curves of $F_N$ is exactly the number of one holed triangles used to construct it. We can compute this number $T_N$ using the Euler characteristic. We begin by observing that $T_N$ must be even because the sides of the one holed triangles are pasted in pairs. There are $N+1$ vertices and $\frac{N(N+1)}{2}$ edges so
$$
N+1 - \frac{N(N+1)}{2} + T_N = 2 - 2 g_{N+1}.
$$
From this 
$$
T_N = 1 - 2 \floor{ \frac{(N-2)(N-3)}{12}} + \frac{N^2}{2}-\frac{N}{2}
$$
which means that $T_N$ grows like $\frac{1}{3} N^2$ in function of $N$. The smallest genus closed surface containing $F_N$ is obtained by pasting the boundary curves in pairs. As $T_N$ is even, this is equal to $g_{N+1} + \frac{T_N}{2}$ so by the above formula is equal to
$$
\frac{N^2}{4}-\frac{N}{2} + \frac{1}{2}. 
$$ 
By a simple topological argument, instead of constructing a minimal genus surface, we can construct a hyperbolic surface, containing $F_N$, of {\it any} genus greater than this minimal genus. To do this, we just attach a hyperbolic surface with a single boundary curve of length $t$ of the appropriate genus to one of the boundary curves and then complete the surface as above. 

To synthesize the above construction, for any $N+1 \equiv 0 \mod 12$ and any $k\geq 0$, there exists a surface with a geometric embedding of $K_{N+1}$ and genus $\frac{N^2}{4} -\frac{N}{2} + \frac{1}{2} + k.$

Now, for any $g\geq 2$, we want to construct a surface of genus $g$ with a geometric embedding of a $K_{N+1}$ for $N$ as large as possible. We choose the maximal $N\equiv -1 \mod 12$ such that 
$$
g \geq \frac{N^2}{2} -\frac{N}{2} + \frac{1}{2} + k
$$
for some $k\geq 0$. Using the condition on the $N$ we are allowed in our construction, we are guaranteed to find a suitable $N$ satisfying
$$
g< \frac{(N+12)^2}{2}-\frac{N+12}{2} + \frac{1}{2} \leq \frac{(N+12)^2}{2} .
$$
From this, we have
$$
N+1\geq \sqrt{ 2 g } - 10
$$
which proves the lower bound in Theorem \ref{thm:B2}.

{\it \underline{Hyperbolic surfaces with infinite chromatic number}}

The building blocks $F_N$ can be used for another purpose - to construct a surface $Z$ such that 
$$
\limsup_{d \to \infty} \chi(Z,d) = \infty.
$$
We note that such a surface must necessarily be of infinite area and in fact the surfaces we shall describe are infinite genus as well. It is entirely possible that there be a much simpler surface with this property, namely $\Hyp$, but this is currently unknown. 

The only thing we require of our surface $Z$ is that it contain copies of $F_{N_k}$ for $N_k \to \infty$ as $k\to \infty$. This is easy to construct as the boundary curves of the $F_N$ all have the same length. One way to do this is to string together the sequence of  $F_{N_k}$, for instance joining one of the boundary geodesics of $F_{N_k}$ with $F_{N_{k+1}}$ for all $k$ (each of them has at least two boundary components so this is possible). We then paste together the remaining (infinite number) of boundary curves in any way. As each surface $F_{N_k}$ has chromatic number $N_k$ for some value of $d$, this proves the result.

\addcontentsline{toc}{section}{References}
\bibliographystyle{amsplain}
\providecommand{\bysame}{\leavevmode\hbox to3em{\hrulefill}\thinspace}
\providecommand{\MR}{\relax\ifhmode\unskip\space\fi MR }
\providecommand{\MRhref}[2]{%
  \href{http://www.ams.org/mathscinet-getitem?mr=#1}{#2}
}
\providecommand{\href}[2]{#2}

\end{document}